\pdfoutput=1 
\documentclass[a4paper,12pt]{article}       
\usepackage{amsmath,amsfonts,amssymb,url,color,epsfig,graphics,float}
\usepackage[justification=centering]{subcaption}
\usepackage{mathrsfs} 
\usepackage{hyperref}
\usepackage[latin1]{inputenc}
\usepackage{anysize}

\newcommand{\R}{{\mathbb{R}}}

\newcommand{\C}{{\mathbb{C}}}

\newcommand{\N}{{\mathbb{N}}}

\def\ha{\frac{1}{2}}

\def\ra{\rightarrow}
\def\ga{\alpha}

\def\gg{\gamma}

\def\gl{\lambda}

\def\gs{\sigma}

\newtheorem{defi}{Definition}[section]
\newtheorem{lemm}[defi]{Lemma}
\newtheorem{prop}[defi]{Proposition}

\newtheorem{coro}[defi]{``Corollary''}
\newtheorem{theo}[defi]{Theorem}

\begin{document}
\title{Propagation of well-prepared states along Martinet singular geodesics}

\author{Yves Colin de Verdi\`ere\footnote{Universit\'e Grenoble-Alpes, 
Institut Fourier, UMR CNRS-UGA 5582, BP74, 38402-Saint Martin d'H\`eres 
Cedex, France (\texttt{yves.colin-de-verdiere@univ-grenoble-alpes.fr})} \ and 
Cyril Letrouit\footnote{Sorbonne Universit\'e, Universit\'e Paris-Diderot, CNRS, 
Inria, Laboratoire Jacques-Louis Lions, F-75005 Paris (\texttt{letrouit@ljll.math.upmc.fr})}\
 \footnote{DMA, \'Ecole normale sup\'erieure, CNRS, PSL Research University, 75005 Paris}}\maketitle
\abstract{We prove that for the Martinet wave equation with ``flat'' metric, which is a subelliptic wave equation,
 singularities can propagate at any speed between $0$ and $1$ along any
  singular geodesic.
 This is in strong contrast with the usual propagation of singularities 
 at speed $1$ for wave equations with elliptic Laplacian.}

\section{Introduction}
\subsection{Propagation of singularities and singular curves}
The celebrated propagation of singularities theorem
 describes the wave-front set $WF(u)$ of a distributional solution $u$ to 
 a partial  differential equation $Pu=f$ in terms of the principal
  symbol $p$ of $P$: it says that, if $p$ is real, then 
  $WF(u)\setminus WF(f)\subset p^{-1}(0)$, and $WF(u)\setminus WF(f)$ is invariant under the 
   bicharacteristic flow induced by the Hamiltonian vector field of $p$.

\smallskip

This result was first proved in \cite[Theorem 6.1.1]{DH} and 
\cite[Proposition 3.5.1]{Ens}. However, it leaves open the case
 where the characteristics of $P$ are not simple, i.e., when there are some points at which $p=dp=0$.
  In a short and impressive paper \cite{Mel86}, Melrose sketched the proof of an analogous 
 propagation of singularities result for the wave operator $P=D_t^2-A$ 
 when $A$ is a self-adjoint non-negative real second-order differential operator
  which is only subelliptic. Such operators $P$ are typical examples for
   which there exist double characteristic points.
   
   \smallskip

Restated in the language of sub-Riemannian geometry (see \cite{Propag}), 
Melrose's result asserts that singularities of subelliptic wave equations propagate
 only along usual null-bicharacteristics (at speed $1$) and along singular
curves (see Definition \ref{d:charac}). Along singular curves, Melrose
writes in  \cite{Mel86}  that the speed should
 be between $0$ and $1$, but nothing more. It is our purpose here to
  prove that for the Martinet wave equation, which is a subelliptic wave equation, singularities 
  can propagate at any speed between $0$ and $1$ along 
  the singular curves of the Martinet distribution. As explained in Point \ref{r:quasi} of Section \ref{s:adaptations}, 
 an analogous result also holds in the so-called quasi-contact case (the computations are
   easier in that case).
  
  \medskip

To state our main result, we consider the \textit{Martinet sub-Laplacian}
$$
\Delta=X_1^2+X_2^2
$$ 
on $\R^3$, where 
 $$
 X_1=\partial_x, \qquad X_2=\partial_y+x^2\partial_z.
 $$
H\"ormander's theorem implies that $\Delta$ is hypoelliptic since $X_1, X_2$ and 
 $[X_1,[X_1,X_2]]$ span $T\R^3$.
  The Martinet
  half-wave equation is
\begin{equation}\label{e:Martinetwave}
i\partial_tu-\sqrt{-\Delta}u=0
\end{equation}
on $\R_t \times\R^3$, with initial datum $u(t=0)=u_0$. The vector fields $X_1$ and $X_2$ span the horizontal
distribution 
$$
\mathcal{D}=\text{Span}(X_1,X_2)\subset T\R^3.
$$
 
Let us recall the definition of singular curves. We use the notation $\mathcal{D}^\perp$ 
for the annihilator   of $\mathcal{D}$ (thus a subcone of the cotangent bundle $T^*\R^3$),  and 
 $\overline{\omega}$ denotes the restriction to $\mathcal{D}^\perp$ of the canonical
  symplectic form $\omega$ on $T^*\R^3$.
\begin{defi} \label{d:charac}
A characteristic curve for $\mathcal{D}$ is an absolutely continuous curve 
$t \mapsto \lambda(t)\in\mathcal{D}^\perp$ that never intersects the zero section of 
$\mathcal{D}^\perp$ and that satisfies $$\dot{\lambda}(t)\in\ker
(\overline{\omega}(\lambda(t)))$$ for almost every $t$. The projection of $\lambda(t)$ onto $\R^3$, which is an 
horizontal curve\footnote{i.e., $d\pi(\dot{\lambda}(t))\in\mathcal{D}_{\lambda(t)}$ for almost every $t$, where $\pi:T^*\R^3\rightarrow \R^3$ denotes the canonical projection.} for $\mathcal{D}$, is called a singular
 curve, and the corresponding characteristic an abnormal extremal lift of that curve.
\end{defi}
We refer the reader to \cite{Mon02} for more material related to sub-Riemannian geometry.

\smallskip

The curve  $t\mapsto \gg(t)=(0,t,0)\in \R^3$ is a singular
curve of the  Martinet distribution $\mathcal{D}$.
 Denoting by $(\xi, \eta, \zeta)$ the dual coordinates  of $(x,y,z)$, 
this curve admits both an abnormal extremal lift, for which $\xi(t)=\eta(t)=0$, and
 a normal extremal lift, for which $\xi(t)=0$, $\eta(t)=1$, $\zeta(t)=0$ (meaning that, if $\tau=1$ is the dual variable 
 of $t$, this  yields a null-bicharacteristic). 
 Martinet-type distributions attracted a lot of attention since Montgomery showed
   in \cite{Mon94} that they provide examples of singular curves which are
   geodesics of the associated sub-Riemannian structure, but which are not necessarily projections of bicharacteristics
   (in contrast with the Riemannian case, where all geodesics are obtained as projections of bicharacteristics).

\smallskip

In this note, all phenomena and computations are done (microlocally) near the 
abnormal extremal lift, and thus away (in the cotangent bundle $T^*\R^3$) from the normal extremal lift,
 which plays no role.

\subsection{Main result}

Let $Y\in C^\infty(\R,\R)$ be equal to $0$ on $(-\infty,1)$ and equal to $1$ on $(2,\infty)$.
Take $\phi \in C_0^\infty (\R,\R)$ with $\phi\geq0$ and $\phi\not\equiv 0$. 
Consider as Cauchy datum for the Martinet half-wave
equation \eqref{e:Martinetwave} the distribution $u_0(x,y,z)$  whose Fourier transform
\footnote{We take the convention $\mathcal{F}f(p)=(2\pi)^{-d}\int_{\R^d}f(q)e^{-iqp}dq$ for the Fourier
 transform in $\R^d$.} with respect to $(y,z)$ is 
\begin{equation} \label{e:Fourier}
 \mathcal{F}_{y,z}u_0(x,\eta,\zeta)=Y(\zeta ) \phi ( \eta /  \zeta^{1/3}) \psi_{\eta,\zeta }(x).
 \end{equation}
Here, $\psi_{\eta,\zeta}$ is the ground state of the $x-$operator
\begin{equation*}
-d_x^2+(\eta+x^2\zeta)^2
\end{equation*}
with $ \psi _{\eta,\zeta}(0)>0$ and $\|\psi_{\eta,\zeta}\|_{L^2}=1$, and  $\alpha_1$ is  the associated eigenvalue. 
Thanks to the Fourier inversion formula applied to \eqref{e:Fourier}, we note that
\[\sqrt{- \Delta }u_0 (x,y,z)=\iint_{\R^2}Y(\zeta ) \phi ( \eta /  \zeta^{1/3})\sqrt{\ga_1}(\eta,\zeta )
 \psi_{\eta,\zeta }(x) e^{i(y\eta +z\zeta )}d\eta d\zeta. \]
We call $u_0$ a {\it well-prepared Cauchy datum.} 
It yields a solution of \eqref{e:Martinetwave}, namely
\[ (U(t)u_0)(x,y,z)=   \iint_{\R^2} Y(\zeta ) \phi ( \eta /  \zeta^{1/3})\psi_{\eta,\zeta }(x)
 e^{-it\sqrt{\alpha_1}(\eta,\zeta )}e^{i(y\eta+ z\zeta )}d\eta d\zeta. \] 
 
 \smallskip
 
For $\mu\in \R$, we set $H_\mu = -d_x^2+ (\mu +x^2)^2 $ and we denote by
$\psi_\mu$ its normalized ground state
$$
H_\mu \psi_\mu = \lambda_1(\mu) \psi_\mu,
$$
whose properties are described at the beginning of Section \ref{sec:quartic}. We also define
$$
F(\mu)=\sqrt{\lambda_1(\mu)}.
$$
We assume that
\begin{equation}\label{e:assumpsizesupport}
F' \text{ is strictly monotonic on the support of $\phi$,}
\end{equation}
 which is no big restriction (choosing adequately the support of $\phi$) since $F$ is an analytic, 
non-affine, function\footnote{See Point 1 of Lemma \ref{lemm:schwartz} and Proposition \ref{p:F}}.
 
 \smallskip
 
We set $\eta= \zeta^{1/3} \eta_1$ and we note that $\psi_{\eta,\zeta}(x)=\zeta^{1/6}\psi_{\eta/\zeta^{1/3}}(\zeta^{1/3}x)=
\zeta^{1/6}\psi_{\eta_1}(\zeta^{1/3}x)$ and $\sqrt{\alpha_1}=\zeta^{1/3}F(\eta/\zeta^{1/3})$. Hence, 
\begin{equation} \label{e:Utu0}
(U(t)u_0)(x,y,z)=   \iint_{\R^2} Y(\zeta ) \zeta^{1/2} \phi ( \eta_1)\psi_{\eta_1 }(\zeta^{1/3}x) 
e^{-i \zeta^{1/3}(tF(\eta _1) -y\eta_1) }e^{i z\zeta }
d\eta_1 d\zeta. 
\end{equation}

We denote by $WF(f)\subset T^*\R^3\setminus 0$ the wave-front set of $f\in\mathcal{D}'(\R^3)$,
 whose projection onto $\R^3$ is the singular support $\text{Sing Supp}(f)$ (see \cite[Definition 8.1.2]{Hor07}).
Our main result states that the speed of propagation of the singularities of $u_0$ is in some window determined by the support of $\phi$.
\begin{theo} \label{t:main} For any $t\in\R$, we have
\begin{equation}\label{e:mainWF} 
WF(U(t)u_0)= \{(0,y,0; \ 0,0,\lambda)\in T^*\R^3, \ \lambda>0, \ y\in tF'(I) \}~,
\end{equation}
where $I$ is the support of $\phi$. In particular,
\begin{equation} \label{e:mainsingsupp}
\text{Sing Supp}(U(t)u_0)= \{(0,y,0)\in\R^3, \ \ y\in tF'(I) \}.
\end{equation}
\end{theo}

Theorem \ref{t:main} means that 
\begin{equation} \label{e:leitmotiv}
\begin{split}
&\textit{singularities propagate   along the singular curve $\gamma$}\\
& \qquad \qquad\textit{ at speeds given by $F'(I)$.}
\end{split}
\end{equation}
Let us comment on the notion of ``speed'' used throughout this paper. In the Riemannian setting, 
when one says that singularities propagate at speed $1$, this has to be understood with respect to the Riemannian metric. 
In the context of the Martinet distribution $\mathcal{D}$, there is also a metric, called sub-Riemannian metric, defined by
\begin{equation} \label{e:metric}
g_q(v)=\inf\left\{u_1^2+u_2^2, \  \ v=u_1X_1(q)+u_2X_2(q)\right\}, \quad q\in\R^3,\ v\in T_q\R^3,
\end{equation}
which is a Riemannian metric on $\mathcal{D}$. This metric $g$ induces naturally a way to measure the speed
  of a point moving along an horizontal curve: if $\delta:J\rightarrow \R^3$ is an horizontal 
  curve describing the time-evolution of a point, i.e., $\dot{\delta}(t)\in \mathcal{D}_{\delta(t)}$ for any $t\in J$,
   then the speed of the point is $(g_{\delta(t)}(\dot{\delta}(t)))^{1/2}$. In the case of the curve $\gamma$, 
   since $g_q(\partial_y)=1$ for any $q$ of the form
 $(0,y,0)$, we have $(g_q(F'(I)\partial_y))^{1/2}=F'(I)$. This is why the set $F'(I)$ is understood as 
 a set of speeds in \eqref{e:leitmotiv}.
 
\begin{prop} \label{p:F'}
There holds $F'(\R)=[a,1)$ for some $-1< a<0$. 
\end{prop}
Together with \eqref{e:leitmotiv}, and choosing $I$ adequately, this implies the following informal statement.
\begin{coro}
 Any value between $0$ and $1$  can be realized as a speed of propagation of singularities along the singular curve $\gamma$.
 \end{coro}
According to \eqref{e:mainWF}, the negative values in the range of $F'$ yield singularities 
propagating backwards along the singular curve. This happens when $F'(I)$ contains negative values (see Proposition \ref{p:F'}).
 
 \paragraph{Organization of the paper.} The paper is organized as follows.
 In Section \ref{s:comm}, we explain in more details the geometrical meaning of  the statement of Theorem \ref{t:main}, 
  and we give several possible adaptations of this result. In Section \ref{sec:quartic}, we prove some properties 
of the eigenfunctions $\psi_\mu$ which play a central role in the next sections. In Section 
\ref{sec:WF0}, we compute the wave-front set of the Cauchy datum $u_0$ thanks to stationary
 phase arguments; this proves Theorem \ref{t:main} at time $t=0$. In Section \ref{sec:WFt}, 
 we complete the proof of Theorem \ref{t:main} by extending the previous computation to any 
 $t\in\R$. We could have directly done the proof for any $t\in\R$ (thus avoiding to distinguish the case $t=0$),
  but we have chosen this presentation to improve readability.
 In Section \ref{sec:functionF}, to illustrate Theorem \ref{t:main}, we prove Proposition \ref{p:F'}, we provide plots of $F$ and $F'$ and compute their asymptotics.

\paragraph{Acknowledgments.} We thank Bernard Helffer and Nicolas Lerner for their help concerning Lemma \ref{lemm:schwartz}. We also thank Emmanuel Tr\'elat for carefully reading a preliminary version of this paper. Finally, we are grateful to an anonymous referee for his questions and suggestions.
 
\section{Comments on the main result}\label{s:comm}
\subsection{Possible adaptations of Theorem \ref{t:main}}\label{s:adaptations}
We describe several possible adaptations of the statement of Theorem \ref{t:main}:
\begin{enumerate}
\item Putting in the initial Fourier data \eqref{e:Fourier} an additional phase $e^{-iz_0\zeta}$ for some fixed $z_0\in\R$,
 we obtain that the singularities of the corresponding solution propagate along the curve $t\mapsto (0,t,z_0)$, which is also
  a singular curve: for this new initial datum, we replace in \eqref{e:mainWF} the $0$ in the $z$ coordinate  by $z_0$.
\item If we consider $(u,D_tu)_{|t=0}=(u_0,0)$ as initial data of the Martinet wave equation $\partial_t^2u-\Delta u=0$,
 the solution is given by 
\[ u(t)=\ha \left( U(t)u_0+ U(-t) u_0 \right). \]
 Hence, under the assumption that $F'(I)$ and $-F'(I)$ do not intersect, \eqref{e:mainWF} must be replaced by
$$
 WF(u(t))=\{(0,y,0; \ 0,0,\lambda)\in T^*\R^3, \ \lambda >0 , \ y\in \pm tF'(I) \}.
$$
\item 
If we take $\zeta<0$ instead of $\zeta>0$ in the (Fourier) initial data $$Y(|\zeta|)
 \phi ( \eta / |\zeta|^{1/3}) \psi_{\eta,\zeta }(x),$$ then we must replace $F'(I)$ by $-F'(-I)$ in the Theorem \ref{t:main}.
The same if we replace $X_2$ by $\partial_y-x^2\partial_z$ and keep $\zeta>0$ in the Fourier 
 initial data. This is due to the ``orientation'' of the singular curve $\gamma$: for Theorem \ref{t:main} to hold without any change,
 we have to take $(0,0,\zeta)(X_2)>0$.
\item \label{r:higheig}
Instead of $\psi_{\eta,\zeta}$, we can use in the Fourier initial datum \eqref{e:Fourier} the $k$-th eigenfunction of
 $-d_x^2+(\eta+x^2\zeta)^2$. This yields a function $F_k$ and the associated velocity $F_k'$, instead of $F$ and $F'$.
 Theorem \ref{t:main} also holds for this initial datum with the same proof, just replacing  $F'$ by $F_k'$ in the statement.
\item \label{r:quasi}
It is possible to establish an analogue of Theorem \ref{t:main} for the half-wave equation
 associated to the quasi-contact sub-Laplacian 
 $$\Delta=\partial_x^2+\partial_y^2+(\partial_z-x\partial_s)^2$$
 on $\R^4$. 
For that, we take Fourier initial data of the form 
 $$\mathcal{F}_{y,z,s}u_0(x,\eta,\zeta,\sigma)= \phi ( \eta /  \gs^{1/2}, 
 \zeta /\gs^{1/2}) \psi_{\eta,\zeta ,\sigma  }(x)$$
where $\phi \in C_0^\infty (\R^2,\R)$, $\eta,\zeta,\sigma$ denote the dual variables of
 $y,z,s$, and $\psi_{\eta,\zeta, \sigma }$ is the normalized ground state of the $x-$operator 
 $-d_x^2+\eta ^2 + (\zeta -x\sigma )^2.$ Then, the singularities propagate along the curve 
 $t\mapsto (0,t,0,0)$ which is a 
 singular curve of the quasi-contact distribution $\text{Span}(\partial_x,\partial_y,\partial_z-x\partial_s)$.
 The proof of this fact requires simpler computations than in the Martinet case since, 
 instead of \emph{quartic} oscillators, they involve usual \emph{harmonic} oscillators. Note that the (non-flat) quasi-contact 
 case has also been investigated in \cite{Sav19}, with other methods.
\end{enumerate}

 \subsection{Geometric comments}

\paragraph{Motivations.}   The result \cite[Theorem 1.8]{Mel86} 
already mentioned in the introduction (and revisited in \cite{Propag})
  implies that in the absence of singular curves, singularities of solutions of the wave equation
  only propagate along null-bicharacteristics. It is in particular the case when the sub-Laplacian has an 
    associated distribution of contact type, since the orthogonal of the distribution
     is in this case symplectic (see \cite[Section 5.2.1]{Mon02}). Another reference for the contact case is \cite{Mel84}.
     Our paper arose from the following questions: 
in the presence of singular curves, can singularities effectively propagate along these singular curves? 
If yes, at which speed(s)?
    
    \smallskip
    
Together with the quasi-contact distribution mentioned
 in Point \ref{r:quasi} of Section \ref{s:adaptations}, 
the Martinet distribution is one of the simplest distributions
 to exhibit singular curves, and this is why we did our computations in this setting.
   
   \smallskip

We now explain that the presence of singular curves for rank $2$ distributions in 3D manifolds is generic. First,
it follows from  Definition \ref{d:charac} that  the existence of singular curves is a property of the distribution $\mathcal{D}$, 
  and does not depend on the metric  $g$ on $\mathcal{D}$ (or on the vector fields $X_1,X_2$ which span $\mathcal{D}$).
Besides, it was proved in \cite[Section II.6]{Mar70} that generically, a rank $2$ distribution $\mathcal{D}_0$
   in a $3D$ manifold $M_0$ is of contact type outside a surface $\mathscr{S}$, called the Martinet surface,
   and near any point of $\mathscr{S}$ except a finite number of them, the distribution is isomorphic to 
   $\mathcal{D}=\text{ker}(dz-x^2dy)$, which is exactly the distribution 
   under study in the present work. Therefore, we expect to be able to generalize Theorem \ref{t:main} 
   to more generic situations of rank $2$ distributions in 3D manifolds.
 
\paragraph{Singular curves as geodesics.}
To explain further the importance of singular curves, let us provide more context about sub-Riemannian geometry.
 A sub-Riemannian manifold is a triple $(M,\mathcal{D},g)$ where $M$ is a smooth 
  manifold, $\mathcal{D}$ is a smooth sub-bundle of $TM$  which is assumed to
  satisfy the H\"ormander condition $\text{Lie}(\mathcal{D}) = TM$, and $g$ is a Riemannian
   metric on $\mathcal{D}$ (which naturally induces a  distance $d$ on $M$). Sub-Riemannian 
   manifolds are thus a generalization of Riemannian manifolds (for which $\mathcal{D}=TM$), 
   and they have been studied in depth since the years 1980, see \cite{Mon02} and \cite{ABB19} for surveys.
   
    \smallskip 
 
Singular curves arise as possible geodesics for the sub-Riemannian distance, i.e.  absolutely continuous
 horizontal paths for which every sufficiently short subarc realizes the sub-Riemannian distance between its endpoints. 
 Indeed, it follows from Pontryagin's maximum principle (see also \cite[Section 5.3.3]{Mon02}) that any sub-Riemannian geodesic is 
 \begin{itemize}
 \item   either \emph{normal}, meaning that it is the projection of an integral curve of the normal  Hamiltonian vector field
 \footnote{By this, we mean the Hamiltonian vector field of $g^*$, the semipositive quadratic form on $T_q^*M$ 
 defined by $g^*(q,p)=\|p_{|\mathcal{D}_q}\|_q^2$, where the norm $\|\cdot \|_q$ is the norm on $\mathcal{D}_q^*$ dual of the norm $g_q$.};
 \item  or \emph{singular}, meaning that it is the projection of a characteristic curve (see Definition \ref{d:charac}).
 \end{itemize}
A sub-Riemannian geodesic can be normal and singular at the same time, and it is indeed the case of the
 singular curve $t\mapsto (x,y,z)=(0,t,0)$ in the Martinet distribution described above. But it was proved in \cite{Mon94}
  that there also exist sub-Riemannian manifolds which exhibit geodesics which are singular, but not normal 
  (they are called strictly singular).

 \smallskip
 
We insist on the fact that in the present work, 
  \begin{center}
\textit{  the minimizing character of the singular curve $\gamma$ plays no role,}
  \end{center}
since our computations can be adapted to the quasi-contact case (see Point \ref{r:quasi} of Section \ref{s:adaptations}), 
 where singular curves are not minimizing.

\paragraph{Spectral effects of singular curves.}
 The study of the spectral consequences of the presence of  singular minimizers was initiated 
in \cite{Mon95}, where it was proved that in the situation where strictly singular minimizers show up as zero loci of 
 two-dimensional magnetic fields, the ground state of a quantum particle concentrates on this curve as $e/h$ tends
  to infinity, where $e$ is the charge and $h$ is the Planck constant. In \cite{CHT-21?}, it is proved that, for 3D compact sub-Riemannian manifolds with 
Martinet singularities, the support of the Weyl measure is the 2D Martinet manifold: most eigenfunctions concentrate on it. 

\smallskip
   
   The present work gives a new illustration of the intuition that singular curves play
   a role ``at the quantum level'', this time at the level of
    propagation for a wave equation. However, the fact that the propagation speed is not $1$, but can take any
     value between $0$ and $1$ was unexpected, since it is in strong contrast with the usual propagation of singularities 
 at speed $1$ for wave equations with elliptic Laplacians.

\section{Some properties of the eigenfunctions $\psi_\mu$}\label{sec:quartic}
Let us recall that $H_\mu $ is the essentially self-adjoint operator\footnote{The operator $H_\mu$ has already been studied for example in \cite{Mon95} and \cite{HP10}.}
$H_\mu =-d_x^2 +(\mu +x^2)^2 $ on $L^2(\R,dx)$ and 
$\psi_\mu $ is the ground state eigenfunction with $\int_\R \psi_\mu(x)^2 dx =1$
and $\psi_\mu (0)>0$. We denote by $\lambda_1(\mu)$ the associated eigenvalue, $\lambda_1(\mu)=F(\mu)^2$. 

\begin{lemm} \label{lemm:schwartz}
The domain of the essentially self-adjoint operator $H_\mu$ is independent of $\mu$. It is denoted by $D(H_0)$. Moreover, the following assertions hold:
\begin{enumerate}
\item The map $\mu \mapsto \gl_1(\mu )$  is analytic on $\R$, and the map $\mu\mapsto \psi_\mu$ is analytic from $\R$ to $D(H_0)$;
\item The function $\psi_\mu $ is in the Schwartz space ${\cal S}(\R)$ uniformly with respect to $\mu $ on any compact subset of $\R$ \footnote{This means that for any compact $K\subset \R$, in the definition of $\mathcal{S}(\R)$, the constants in the semi-norms can be taken independent of $\mu\in K$.};
\item Any derivative in $D(H_0)$ of the map $\mu \mapsto \psi_ \mu$ is in the Schwartz space ${\cal S}(\R)$ uniformly with respect to $\mu $ on any compact subset of $\R$.
\end{enumerate}
\end{lemm}
\textit{Proof.}
 The domain of $H_\mu$ is given by 
\begin{align*}
 D(H_\mu )&=\{ \psi \in L^2(\R), -\psi'' +(\mu+x^2)^2\psi  \in L^2(\R)\} \\
&=\{ \psi \in L^2(\R), -\psi'' +x^4 \psi  \in L^2(\R),\ x^2 \psi \in L^2(\R)\} \\
&= \{\psi \in L^2(\mathbb{R}), - \psi'' \in L^2(\mathbb{R}),\ x^4 \psi \in L^2(\mathbb{R}),\ x^2 \psi \in L^2(\mathbb{R})\}
 \end{align*}
(see \cite{Sim70} and \cite{EG74}).
We have hence $D(H_\mu )=D(H_0)$. The map  $\mu \mapsto H_\mu $ is analytic from $\R$ into
$\mathcal{L}(D(H_0), L^2(\R))$. Moreover, by \cite[Theorem 3.1]{BS12}, the eigenvalues of $H_\mu $ are non-degenerate 
(simple). This implies (see \cite[Chapter VII.2]{Kat13} or \cite[Proposition 5.25]{CR19}) that  the eigenvalues $\gl_1(\mu )$  and eigenfunctions $\psi_\mu $
  are analytic functions of $\mu $, respectively 
with values in $\R$ and in $D(H_0 )$. This proves Point 1.

Point 2 follows from Agmon estimates (precisely, \cite[Proposition 3.3.4]{Hel88} with $h=h_0=1$), which are uniform with respect to $\mu$ on any compact subset of $\R$.

This allows to start to prove Point 3 by induction. Assume that Point 3 is true for the derivatives of order $0,\ldots,k-1$. Then,
taking the derivatives with values in the domain $D(H_0)$ with respect to $\mu $ in the equation $(H_\mu-\lambda_1(\mu))\psi_\mu=0$, we get
\begin{equation} \label{e:res} (H_\mu -\gl_1(\mu) ) \frac{d^k}{d\mu ^k}\psi_\mu = v_{k,\mu} \end{equation}
and we know, by the induction hypothesis, that $v_{k,\mu} \in {\cal S}(\R)$ uniformly with respect to $\mu$ on any compact subset of $\R$. 
We now use the results of \cite[Section 25]{Shu87} (see also \cite[Section 23]{Shu87} for the notations, and \cite{HR82} for similar results).
We check that $\xi^2 + x^4$ is a symbol in the sense of Definition 25.1 of \cite{Shu87}, with $m=4$, $m_0=2$ and $\rho=1/2$.
Its standard quantization (i.e., $\tau=0$ in Equation (23.31) of \cite{Shu87}) is $H_\mu$.
 By \cite[Theorem 25.1]{Shu87}, $H_\mu-\lambda_1(\mu)$ admits a parametrix $B_\mu$; 
in particular, $B_\mu(H_\mu-\lambda_1(\mu))=\text{Id}+R_\mu$  where $R_\mu$ is smoothing. Hence, composing on the left by $B_\mu$ in \eqref{e:res}, 
and noting that $B_\mu v_{k,\mu}\in {\cal S}(\R)$, we obtain that $\frac{d^k}{d\mu ^k}\psi_\mu\in{\cal S}(\R)$ uniformly with respect to $\mu$ 
on any compact subset of $\R$, which concludes the induction and the proof of Point 3.

\section{Wave-front of the Cauchy datum}  \label{sec:WF0}
The goal of this section is to compute the wave-front set of $u_0$. In other words, we prove Theorem \ref{t:main} for $t=0$. Recall that (see \eqref{e:Utu0})
\begin{equation}  \label{e:defu0}
u_0(x,y,z)=\iint_{\R^2} Y(\zeta ) \zeta^{1/2} \phi ( \eta_1) \psi_{\eta_1}(\zeta^{1/3}x) e^{i(y\zeta^{1/3}\eta_1 +z\zeta )}d\eta_1 d\zeta. 
\end{equation}

\begin{lemm}  \label{l:singsupp}
The function $u_0$ is smooth on $\R^3\setminus \{(0,0,0)\}$.
\end{lemm}
\textit{Proof.} We prove successively that $u_0$ is smooth outside $x=0$, $y=0$ and $z=0$. Any derivative of
 \eqref{e:defu0} in $x,y,z$ is of the form
\begin{equation} \label{e:der}  
\iint_{\R^2} Y(\zeta ) \zeta^{\alpha} \psi_{\eta_1}^{(\gamma)}(\zeta^{1/3}x) \phi ( \eta_1)
\eta_1^{\beta}e^{i(y\zeta^{1/3}\eta_1 +z\zeta )}d\eta_1 d\zeta 
\end{equation}
for some $\alpha,\beta,\gamma\geq 0$. By the dominated convergence theorem, locally uniform (in $x,y,z$) 
convergence of these integrals implies smoothness.
 Recalling that $\phi$  has compact support, we see that the main difficulty for proving smoothness comes from the integration
 in $\zeta$ in \eqref{e:der}.
 
 \smallskip

For $x\ne 0 $ it follows from Lemma \ref{lemm:schwartz} (Point 2) that the integrand in \eqref{e:der}
 has a fast decay in $\zeta$. This proves that $u_0$ is smooth outside $x=0$.
 
 \smallskip

If $y\ne 0$, we use the fact that the phase $y\zeta^{1/3}\eta_1+z\zeta $ is non critical 
with respect to $\eta _1$ to get the decay in $\zeta $.  More precisely, \eqref{e:der} is equal to
$$
\iint_{\R^2} Y(\zeta ) \zeta^{\alpha}  (y\zeta^{1/3})^{-N}D_{\eta_1}^N(\psi_{\eta_1}^{(\gamma)}(\zeta^{1/3}x)
\phi ( \eta_1)\eta_1^{\beta})e^{i(y\zeta^{1/3}\eta_1 +z\zeta )}d\eta_1 d\zeta 
$$
after integration by parts in $\eta_1$ (where $D_{\eta_1}=i^{-1}\partial_{\eta_1}$). Taking $N$ sufficiently large and using that
$D_{\eta_1}^N(\psi_{\eta_1}^{(\gamma)}(\zeta^{1/3}x)
\phi ( \eta_1)\eta_1^{\beta})$ is bounded thanks to Lemma \ref{lemm:schwartz} (Point 3), we obtain that this integral 
converges when $y\neq0$, and that this convergence is locally uniform with respect to $x,y,z$.
 This proves that $u_0$ is smooth outside $y=0$.

\smallskip

Finally, let us study the case $z\neq 0$. We can also assume that $y\leq 1$ due to the previous point.

\textit{Claim.} The function 
\begin{equation} \label{e:symbolzeta}\zeta \mapsto Y(\zeta ) \zeta^{1/2} \phi ( \eta_1)\psi_{\eta_1}^{(\gamma)}(\zeta^{1/3}x)
e^{iy\zeta^{1/3}\eta_1} \end{equation}
is a symbol (see Definition \ref{d:symbol}) uniformly on every compact in $(y,\eta _1)$. 

\textit{Proof.} The functions $\zeta\mapsto \zeta^{1/2}\phi ( \eta_1)$ and $\zeta \mapsto Y(\zeta ) e^{iy\zeta^{1/3}\eta_1}$ 
 are symbols (with $\rho=1$ and $\rho=2/3$ respectively, see Definition \ref{d:symbol}) uniformly on every compact in $(y,\eta_1)$. Besides, $\zeta \mapsto \psi_{\eta_1}^{(\gamma)}(\zeta^{1/3}x)$
 is also a symbol (of degree $0$ with $\rho=1$): we notice for example that the first derivative with respect to $\zeta $ writes
$(1/3 )\zeta^{-1}(\zeta^{1/3} x) \psi^{(\gamma+1)}_{\eta _1}(\zeta^{1/3}x)$ which is uniformly $O(1/\zeta )$
 thanks to Lemma \ref{lemm:schwartz} (Point 2). 
Finally, since the space of symbols is an algebra for the pointwise product, we get the claim.

Integrating \eqref{e:symbolzeta} in $\eta_1\in\R$ and using Lemma \ref{lemm:symbol} (in the variable $\zeta$), we obtain 
 that \eqref{e:defu0} is smooth outside $z=0$, which concludes the proof of Lemma \ref{l:singsupp}.

\medskip

The following lemma proves Theorem \ref{t:main} at time $t=0$.

\begin{lemm}  \label{l:WF}
There holds $WF(u_0)=\{(0,0,0; \ 0,0,\lambda)\in T^*\R^3, \ \lambda>0\}$.
\end{lemm}
\textit{Proof.}
The Fourier transform of $u_0$ is 
\begin{equation} \label{e:Fouriertransform} U_0(\xi,\eta,\zeta)=Y(\zeta )
\phi ( \eta /  \zeta^{1/3}) \Psi_{\eta/\zeta^{1/3}}(\xi/\zeta^{1/3}) 
\end{equation}
where $\Psi_\mu $ is the Fourier transform of the eigenfunction $\psi_\mu $. By Lemma \ref{lemm:schwartz} (Point 2), for any $N\in\N$ we get 
\begin{equation} \label{e:decayU0}
|U_0(\xi,\eta,\zeta)|\leq C_N |\phi ( \eta /  \zeta^{1/3})| (1+|\xi/\zeta^{1/3}|)^{-N}. 
\end{equation}

We show that $U_0$ is fastly decaying in any cone $C:=\{ |\xi |+ |\eta |\geq c|\zeta | \}$ for $c$ small.
We split the cone 
into $C=C_1\cup C_2$ with
$C_1=C \cap \{ |\xi |\leq |\eta |\} $
and $C_2=C \cap \{ |\eta |\leq |\xi |\} $.\\
In $C_1$, we have $|\eta/\zeta^{1/3}|\geq c_1|\eta^{2/3}|$. This implies that $\phi(\eta/\zeta^{1/3})$
 vanishes for $\eta$ large enough. Hence, $U_0$ has fast decay in $C_1$.\\
 In $C_2$,  we have $| \xi/\zeta^{1/3}|\geq c_2 |\xi|^{2/3}\geq c_3(1+\xi^2+\eta^2+\zeta^2)^{1/3}$,
 hence, plugging into \eqref{e:decayU0}, we get that $U_0$ has fast decay in $C_2$.

This proves that no point of the form $(x,y,z; \ \xi,\eta,\zeta)\in T^*\R^3$ with $(\xi,\eta)\neq (0,0)$
 can belong to $WF(u_0)$. Moreover, due to the factor $Y(\zeta)$, necessarily $WF(u_0)\subset \{\zeta>0\}$.
 Combining with Lemma \ref{l:singsupp}, we get the inclusion $\subset$ in Lemma \ref{l:WF}.

\smallskip

Let us finally prove that $(0,0,0; \ 0,0,\lambda) \in WF(u_0)$ for $\lambda>0$. We pick $a,b\in\R$ such that $\phi(a)\neq 0$
 and $\Psi_a(b)\neq0$. Then, we note that $U_0(\zeta^{1/3}a,\zeta^{1/3}b,\zeta)$ is independent of $\zeta$
 and $\neq 0$, thus it is not fastly decaying as $\zeta\rightarrow +\infty$. Since $(\zeta^{1/3}a,\zeta^{1/3}b,\zeta)$
 converges to the direction $(0,0,+\infty)$ as $\zeta\rightarrow +\infty$, we get that there exists at least
 one point of the form $(x,y,z; \ 0,0,\lambda)\in T^*\R^3$ which belongs to $WF(u_0)$. 
By Lemma \ref{l:singsupp}, we necessarily have $x=y=z=0$, which concludes the proof.

\section{Wave front of the propagated solution} \label{sec:WFt}
In this Section, we complete the proof of Theorem \ref{t:main}. We set
\begin{equation*}
\mathcal{G}_t=\{(0,y,0; \ 0,0,\lambda), \ \lambda>0, \ y\in tF'({\rm Support }(\phi )) \}.
\end{equation*}
In Section \ref{sec:subset}, we prove the inclusion $WF(U(t)u_0)\subset \mathcal{G}_t$, and then in Section \ref{sec:egal}
 the converse inclusion $\mathcal{G}_t\subset WF(U(t)u_0)$. This completes the proof of Theorem \ref{t:main}.

\subsection{The inclusion $WF(U(t)u_0)\subset \mathcal{G}_t$} \label{sec:subset}

For this inclusion, we follow the same arguments as in Section \ref{sec:WF0}: we adapt Lemma \ref{l:singsupp}
 to find out the singular support of $U(t)u_0$, and then we adapt Lemma \ref{l:WF} to determine the full wave-front set.
\begin{lemm}  \label{l:singsupp2}
For any $t\in\R$, $U(t)u_0$ is smooth outside $\{(0,y,0)\in \R^3,\ y\in tF'(I)\}$.
\end{lemm}
\textit{Proof.} As in Lemma \ref{l:singsupp}, we prove successively that $U(t)u_0$ is smooth outside $x=0$, $y\notin tF'(I)$ and $z=0$.
 Any derivative of $U(t)u_0$ is of the form
\begin{equation} \label{e:der2}
\iint_{\R^2} Y(\zeta ) \zeta^{\alpha} \psi_{\eta_1}^{(\gamma)}(\zeta^{1/3}x)
 \phi ( \eta_1)\eta_1^{\beta}e^{-i\zeta^{1/3}(tF(\eta_1)-y\eta_1)}e^{iz\zeta }d\eta_1 d\zeta
\end{equation}
for some $\alpha,\beta,\gamma\geq 0$.

\smallskip

For $x\ne 0 $, it follows from Lemma \ref{lemm:schwartz} (Point 2) that the integrand in \eqref{e:der2}
 has a fast decay in $\zeta$ (locally uniformly in $x,y,z$). This proves that $U(t)u_0$ is smooth outside $x=0$.

\smallskip

 If $y\notin tF'(I)$, we use the fact that the phase $\zeta^{1/3}(tF(\eta _1) -y\eta_1) - z\zeta$ is non critical with respect to $\eta _1$ to get decay in $\zeta $.
We set $R_{\eta_1} H=D_{\eta_1}(Q^{-1}H)$ where $Q=D_{\eta_1}(-i(\zeta^{1/3}(tF(\eta_1)-y\eta_1)-z\zeta))=-\zeta^{1/3}(tF'(\eta_1)-y)$. Note that $Q\neq 0$ since $y\notin tF'(I)$. Doing $N$ integration by parts, the above expression becomes
\begin{equation} \label{e:IPP2}
\iint_{\R^2} Y(\zeta ) \zeta^{\alpha}  R_{\eta_1}^N(\psi_{\eta_1}^{(\gamma)}(\zeta^{1/3}x)\phi ( \eta_1)\eta_1^{\beta})e^{-i\zeta^{1/3}(tF(\eta_1)-y\eta_1)}e^{iz\zeta }d\eta_1 d\zeta.
\end{equation}
We set $H(x,\eta_1,\zeta)=\psi_{\eta_1}^{(\gamma)}(\zeta^{1/3}x)\phi ( \eta_1)\eta_1^{\beta}$.
\smallskip

\textit{Claim.} For any $N$, there exists $C_N$ such that $|R_{\eta_1}^N H(x,\eta_1,\zeta)|\leq C_N|\zeta|^{-N/3}$ for any $\zeta\in\R$, any $\eta_1\in I=\text{Support}(\phi)$ and any $x\in\R$.

\smallskip

Taking $N$ sufficiently large, the claim implies that \eqref{e:IPP2}, and thus \eqref{e:der2}, 
converge (locally uniformly), which proves the smoothness when $y\notin tF'(I)$
 thanks to the dominated convergence theorem.

\smallskip

\textit{Proof of the claim.} We prove it first for $N=1$. We have
\begin{equation} \label{e:reta1}
R_{\eta_1}H=\frac{D_{\eta_1}H}{Q}-H\frac{D_{\eta_1}Q}{Q^2}.
\end{equation}
Since $H$ is bounded (thanks to Point 2 of Lemma \ref{lemm:schwartz}) and $|Q|\geq c|\zeta|^{1/3}$ with $c>0$ and $|D_{\eta_1}Q|\leq C|\zeta|^{1/3}$ on the support of $\phi$, we have $|H\frac{D_{\eta_1}Q}{Q^2}|\leq c|\zeta|^{-1/3}$. For the first term in the right-hand side of \eqref{e:reta1}, we only need to prove that $D_{\eta_1}H$ is bounded. When $D_{\eta_1}$ falls on $\phi(\eta_1)$ or $\eta_1^\beta$, it is immediate. When $D_{\eta_1}$ falls on $\psi_{\eta_1}^{(\gamma)}(\zeta^{1/3}x)$, we use Lemma \ref{lemm:schwartz} (Point 3) and also get the result. This ends the proof of the case $N=1$. Now, we notice that our argument works not only for $H$,  but for any function of the form $\psi_{\eta_1}^{(\gamma')}(\zeta^{1/3}x)\phi^{(\delta)}(\eta_1)\eta_1^{\beta'}$ where $\phi^{(\delta)}$ is any derivative of $\phi$ and $\beta',\gamma'\geq 0$. Hence, applying the previous argument recursively, we obtain the claim for any $N$.

\medskip

Finally, the case   $z\neq 0$ is checked in the same way as in the case $t=0$, just shifting the phase by $it\zeta^{1/3}F(\eta_1)$ in \eqref{e:symbolzeta}.

\medskip

{\it Let us finish  the proof of the inclusion $WF(U(t)u_0)\subset \mathcal{G}_t$.}

 The Fourier transform of $U(t)u_0$ is 
\begin{equation}\label{equ:ut}
 \mathcal{F}(U(t)u_0)(\xi,\eta,\zeta)=Y(\zeta )
\phi ( \eta /  \zeta^{1/3}) \Psi_{\eta/\zeta^{1/3}}(\xi/\zeta^{1/3})e^{-it\sqrt{\alpha_1}(\eta,\zeta )}. 
\end{equation}
 The change of phase with respect \eqref{e:Fouriertransform} has no influence on the properties of decay at infinity. 
 Hence, the proof of Lemma \ref{l:WF} allows to conclude that $WF(U(t)u_0)\subset \mathcal{G}_t$ for any $t\in\R$. 

\subsection{The inclusion $\mathcal{G}_t\subset WF(U(t)u_0)$} \label{sec:egal}
We fix $t\in\R$ and we prove the non smoothness at $(0,tF'(c) ,0) $ for any $c\in I$.
We can assume that $c$ is in the interior of $I$ and that $\phi(c)\neq 0$. 
This implies thanks to \eqref{e:assumpsizesupport} that $F''(c)\neq 0$.
We want to  show non-smoothness with respect to $z$ at $x=0$,  $y=tF'(c) $ and $z=0$. We set $v(z):=(U(t)u_0)(0,tF'(c),z)$.
We will show that the Fourier transform of $v$ is not fastly decaying.

\smallskip

Starting from \eqref{e:Utu0}, we get the explicit formula for the Fourier transform of $v$,
\begin{equation*}
\mathcal{F}v(\zeta)=  Y(\zeta ) \zeta^{1/2} K(\zeta) 
\end{equation*}
where
\begin{equation*}
K(\zeta)=\int_{\R}  \phi ( \eta_1)\psi_{\eta_1 }(0) 
e^{-i \zeta^{1/3}t(F(\eta _1) -F'(c)\eta_1) } d\eta_1.
\end{equation*}
The only critical point of the phase $\eta_1\mapsto -i \zeta^{1/3}t(F(\eta _1) -F'(c)\eta_1)$ located in $I$ is $c$ thanks to \eqref{e:assumpsizesupport}.
Applying the stationary phase theorem with respect to $\eta_1$, we obtain
\begin{equation*}
K(\zeta)=e^{-i\zeta^{1/3}t(F(c)-F'(c)c)} \sum_{j\geq 1}a_j(\zeta^{1/3}|t|)^{-j/2}
\end{equation*}
where 
$$
a_1=\phi(c)\psi_{c}(0)\left(\frac{2\pi}{|F''(c)|}\right)^{1/2}\exp(-i\frac{\pi}{4}\text{sgn}(F''(c)))\neq 0.
$$
Since $\phi(c)> 0$ and $\psi_{c}(0)>0$,we have $K(\zeta)\sim c_0(\zeta^{1/3}|t|)^{-1/2}$ where $c_0\neq 0$, and $\mathcal{F}v(\zeta)$ is not fastly decaying as $\zeta\rightarrow +\infty$. Applying Lemma \ref{lemm:symbol} to $a=\mathcal{F}v$ which is a symbol in $\zeta$, this implies that $v$ is not smooth at $z=0$, thus $U(t)u_0$ is not smooth at $(0,tF'(c),0)$. 

\section{The function $F_k(\mu)=\sqrt{\gl_k(\mu )}$} \label{sec:functionF}
In this Section, we illustrate Theorem \ref{t:main}
 with some plots and asymptotics of the functions $F_k$ defined by  $\mu \ra \sqrt{\gl_k(\mu )}$. 
As shown by Theorem \ref{t:main} (and Point \ref{r:higheig} in Section \ref{s:adaptations}), the speeds of the propagation of singularities
 along the singular curve are determined by the derivative $F'_k(\mu)$. Below,
 we plot $F=F_1$ and $F'$ for $\mu\in(-10,10)$\footnote{We thank Julien Guillod for his help in making the first numerical experiments.}. 
\begin{figure}[H] 
\captionsetup[subfigure]{justification=centering}
    \begin{subfigure}[b]{0.3\textwidth}
        \subcaptionbox{Plot of $F(\mu)$ for $\mu\in (-10,10)$}{ \includegraphics[width=8cm]{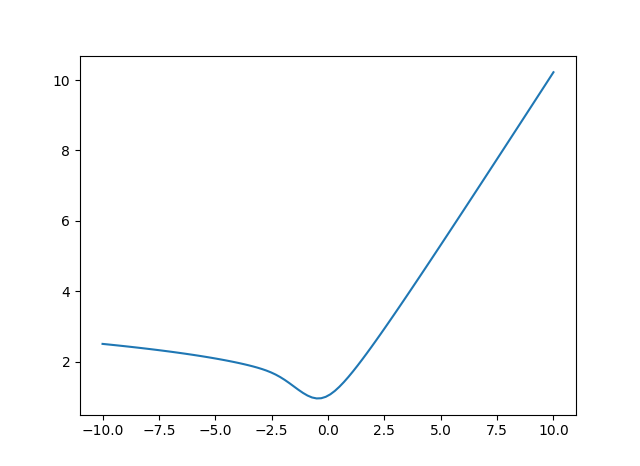}}
        \label{fig:F}
    \end{subfigure}
\qquad \qquad \qquad \qquad
    \begin{subfigure}[b]{0.3\textwidth}
       \subcaptionbox{Plot of $F'(\mu)$ for $\mu\in (-10,10)$}{ \includegraphics[width=8cm]{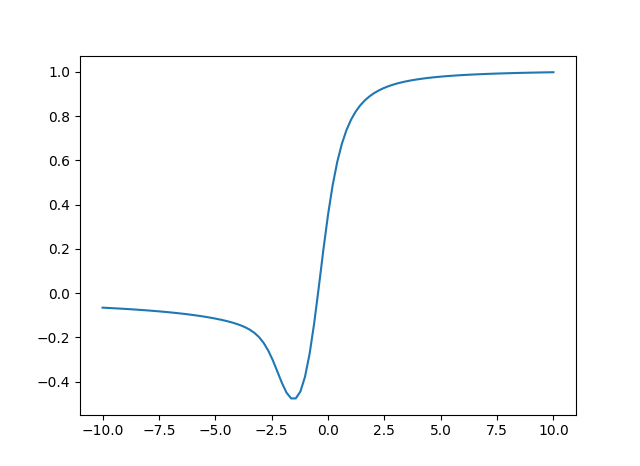}}
        \label{fig:F'}
    \end{subfigure}
\end{figure}

Recall that the $F_k$'s are  analytic (see Point 1 of Lemma \ref{lemm:schwartz}). We state a more precise version of Proposition \ref{p:F'}:
\begin{prop} \label{p:F}
For any $k\in\N\setminus\{0\}$, there holds $F'_k(\mu)\rightarrow 1^-$ as $\mu\rightarrow +\infty$, $F'_k(\mu)\ra 0^-$
 as $\mu\rightarrow -\infty$, and $F_k'$ is minimal  for some value $\mu_k^\star <0$.
 There exists $a_k\in (-1,0)$ such that the range of $F'_k$ is $[a_k,1)$.
\end{prop}
Proposition \ref{p:F} will be a consequence of the following result:  
\begin{prop} \label{p:helffer}
Denote by $\lambda_k(\mu)$ the $k$-th eigenvalue of $H_\mu=-d_x^2+(\mu+x^2)^2$. Then, for $k\in\N\setminus \{0\}$, as $\mu \ra +\infty $,
\begin{equation}\label{e:lambdak}
 \gl_k(\mu )=\mu^2 + \sqrt{2}(2k-1)\sqrt{\mu }+ \sum_{\ell=2}^\infty  b_{\ell,k}\mu ^{2-3\ell/2} 
 \end{equation}
and 
\begin{equation}\label{e:derlambdak}
 \frac{d}{d\mu }\sqrt{\gl_k(\mu )}= 1-\frac{2k-1}{2\sqrt{2}}\mu ^{-3/2}+O(\mu^{-3})
 \end{equation}
These derivatives are  $>0$ and converge to $1$. 

As $\mu \ra -\infty $, for $k\in\N\setminus \{0\}$, 
\begin{align}
& \gl_{2k-1}(\mu )= 2(2k-1)\sqrt{-\mu }+ \sum_{\ell=2}^\infty  c_{\ell,k}(-\mu )^{2-3\ell/2}  \label{e:lambdaimpair}\\
& \gl_{2k}(\mu )=\gl_{2k-1}+ o\left(\mu^{-\infty }\right) \label{e:lambdapair}
 \end{align}
and
\begin{align}
 \frac{d}{d\mu }\sqrt{\gl_{2k-1}(\mu )}=-\frac{\sqrt{2(2k-1)}}{4}(-\mu)^{-3/4}+O(\mu^{-3/2}) \label{e:derlambdaimpair}
 \end{align}
and the same for $\frac{d}{d\mu }\sqrt{\gl_{2k}(\mu)} $. 
These derivative are  $<0$ and converge to $0$. 
\end{prop}

\textit{Proof of Proposition \ref{p:helffer}.} For $\mu>0$, we consider the operator
 $T_\mu:\psi\mapsto\psi(\cdot/\mu^{1/4})$. Then $H_\mu=T_\mu^{-1} G_\mu T_\mu$ where $G_\mu=\mu^2+\mu^{1/2}(-d_x^2+2x^2+x^4/\mu^{3/2})$.
 The eigenvalues of $-d_x^2+2x^2+hx^4$ for $h\rightarrow 0$ can be computed with the usual perturbation theory
 (see \cite[Chapter XII.3]{RS78}), and this yields \eqref{e:lambdak} with $h=\mu^{-3/2}$. Moreover the formal expansion can be
differentiated with respect to $\mu $, hence we get \eqref{e:derlambdak}.

For $\mu=-{\mu_0}<0$, we see that the transformation $x\mapsto \mu_0^{1/4}(x\mp \mu_0^{1/2})$ conjugates $H_\mu$
 to the operator $\mu_0^{1/2}(-d_x^2+4x^2\pm 4\mu_0^{-3/4}x^3+\mu_0^{-3/2}x^4)$. Using again perturbation theory and
 the separation into pairs of eigenvalues in double wells (see \cite{HS84}), we get \eqref{e:lambdaimpair} and
 \eqref{e:lambdapair}, and \eqref{e:derlambdaimpair} follows. 
\medskip

\textit{Proof of Proposition \ref{p:F}.} The convergences at $\pm \infty$ are proved by Proposition \ref{p:helffer}.
 This behaviour at $\pm \infty$ implies the existence of $\mu_k^\star $ such that $F'_k(\mu_k^\star)=a_k$
is minimal. We denote by $\psi_\mu^k$ the normalized eigenfunction corresponding to $\lambda_k(\mu)$. 
Taking the first derivative (with value in the domain $D(H_0)$) with respect to $\mu $ of the
 eigenfunction equation $(H_\mu-\lambda_k(\mu))\psi_\mu^k=0$, and then integrating against $\psi_\mu^k$, 
 we obtain  $\lambda_k'(\mu)=\int_{\R} (\mu+x^2)\psi_\mu^k(x)^2dx$. 
Thus, $$F'_k(\mu)=\frac{1}{\sqrt{\lambda_k(\mu)}}\int_{\R} (\mu+x^2)\psi_\mu^k(x)^2dx $$ 
which is positive for $\mu \geq 0$, hence $\mu_k^\star<0$. 

\smallskip

It remains to show that $|F'_k(\mu )|<1$ for every $\mu $: by the Cauchy-Schwarz inequality, we get
 $$F'_k(\mu)^2\leq \frac{1}{{\lambda_k(\mu)}}\int_{\R} (\mu+x^2)^2\psi^k_\mu(x)^2dx \int_{\R} \psi^k_\mu(x)^2dx $$ 
and, from the quadratic form associated to $H_\mu$,
$$ \int_{\R} (\mu+x^2)^2\psi_\mu^k(x)^2dx < \lambda_k(\mu),$$ 
which concludes the proof.
\section*{Appendix}

\appendix

\section{Fourier transform of symbols} \label{sec:fourier}

\begin{defi} \label{d:symbol}
A smooth function $a:\R^d \ra \C $ is called a {\it symbol} of degree $\leq m $ if there exists  
 $0< \rho \leq 1$  so that
the partial derivatives of $a$ satisfy
\[  \forall \ga\in\N^d , \qquad |D^\ga a (\xi) |\leq C_\ga (1+|\xi|)^{m-\rho |\ga |}.\]
\end{defi}
The space of symbols is an algebra for the pointwise product. If
$a$ is a real valued symbol of degree $m<1$ and $\rho >m$, 
$e^{ia } $ is a symbol of degree $0$ (with a different $\rho$).

\smallskip

We will need the
\begin{lemm}\label{lemm:symbol}
If $a$ is a symbol, the Fourier transform  $\mathcal{F}a$ of 
 $a$ is smooth outside $x=0$ and all derivatives of $\mathcal{F}a$ decay fastly at  infinity. 
If moreover $a$ does not belong to the Schwartz space ${\cal S}(\R^d)$, then $\mathcal{F}a$ is non smooth at $x=0$.

\end{lemm}
\textit{Proof.} 
For $x\neq 0$ and for any $\alpha, \beta\in\N^d$, we have
\begin{equation} \label{e:hata}
(\mathcal{F}a)^{(\beta)}(x)=C_\beta\int_{\R^d}\xi^\beta a(\xi)e^{-ix\xi}d\xi=\frac{c_\beta^\alpha}{x^\alpha} \int_{\R^d}
 D^\alpha_\xi(\xi^\beta a(\xi))e^{-ix\xi}d\xi.
\end{equation}
The multi-index $\beta\in \N^d$ being fixed, this last integral converges for $|\alpha|$ sufficiently large since $a$ is a symbol. 
By the dominated convergence theorem, this implies that $\mathcal{F}a$ is smooth outside $x=0$.
Moreover, \eqref{e:hata} also implies that all derivatives of $\mathcal{F}a$ decay fastly at infinity.

Finally, if  $\mathcal{F}a$ were  smooth at $0$, then $\mathcal{F}a$ would be in the Schwartz space as well as $a$.

\end{document}